\documentclass[11pt]{article}
\usepackage{amsmath, amsthm, amssymb}    % May not all be necessary
\usepackage{bridges}
\usepackage{graphicx,xcolor}			% For including pictures
\usepackage[colorlinks=true, urlcolor=blue, citecolor=black, linkcolor=black]{hyperref}  
\usepackage[export]{adjustbox}
\usepackage[font={small}]{caption}
\usepackage{subcaption}			% For captioning multi-panel figures

\usepackage[nameinlink,capitalize,noabbrev]{cleveref}
\usepackage{paralist}

\newtheorem{theorem}{Theorem}
\newtheorem*{theorem*}{Theorem}

\newtheorem{proposition}{Proposition}

\theoremstyle{remark}

\theoremstyle{definition}
\newtheorem{definition}{Definition}

\def\cn{\mathrm{cn}}
\def\sn{\mathrm{sn}}

\def\dn{\mathrm{dn}}

\hypersetup{
    pdftoolbar=true,        % show Acrobat’s toolbar?
    pdfmenubar=true,        % show Acrobat’s menu?
    pdffitwindow=false,     % window fit to page when opened
    pdfstartview={FitH},    % fits the width of 
    colorlinks=true,       % false: boxed links; true: colored links
    linkcolor=olive,          % color of internal links (change box color with linkbordercolor)
    citecolor=blue,        % color of links to bibliography
    filecolor=black,      % color of file links
    urlcolor=red           % color of external links
}

\renewcommand{\H}{\mathcal{H}}

\newcommand{\K}{\mathcal{K}}
\renewcommand{\S}{\mathcal{S}}

\newcommand{\E}{\mathcal{E}}

\renewcommand{\H}{\mathcal{H}}

\newcommand{\rulesep}{\unskip\ \vrule\ }
\usepackage{esvect}

\title{Poncelet Spatio-Temporal Surfaces and Tangles}

\author{Claudio Esperança\textsuperscript{1}, Ronaldo Garcia\textsuperscript{2}, and Dan Reznik\textsuperscript{3}
\vspace{10pt}\\
\textsuperscript{1}PESC, Fed. Univ. Rio, Rio de Janeiro, Brazil; claudio.esperanca@gmail.com\\
\textsuperscript{2}Math and Stats dept., Fed. Univ. Goiás, Goiânia, Brazil; ragarcia@ufg.br\\
\textsuperscript{3}Data Science Consulting, Rio de Janeiro, Brazil; dreznik@gmail.com}
\date{}	

\begin{document}

\maketitle
% Prevent page number 1 from being printed on the first page.
%\thispagestyle{empty}

\begin{abstract}
We explore geometric properties of 3d surfaces swept by a family of Poncelet triangles, as well as tangles produced by space curves they define.
\end{abstract}

\section{Introduction}
\label{sec:intro}
\begin{figure}
 \centering
     \begin{subfigure}[c]{0.49\textwidth}
         \centering
         \includegraphics[trim=0 0 0 0,clip,width=\textwidth]{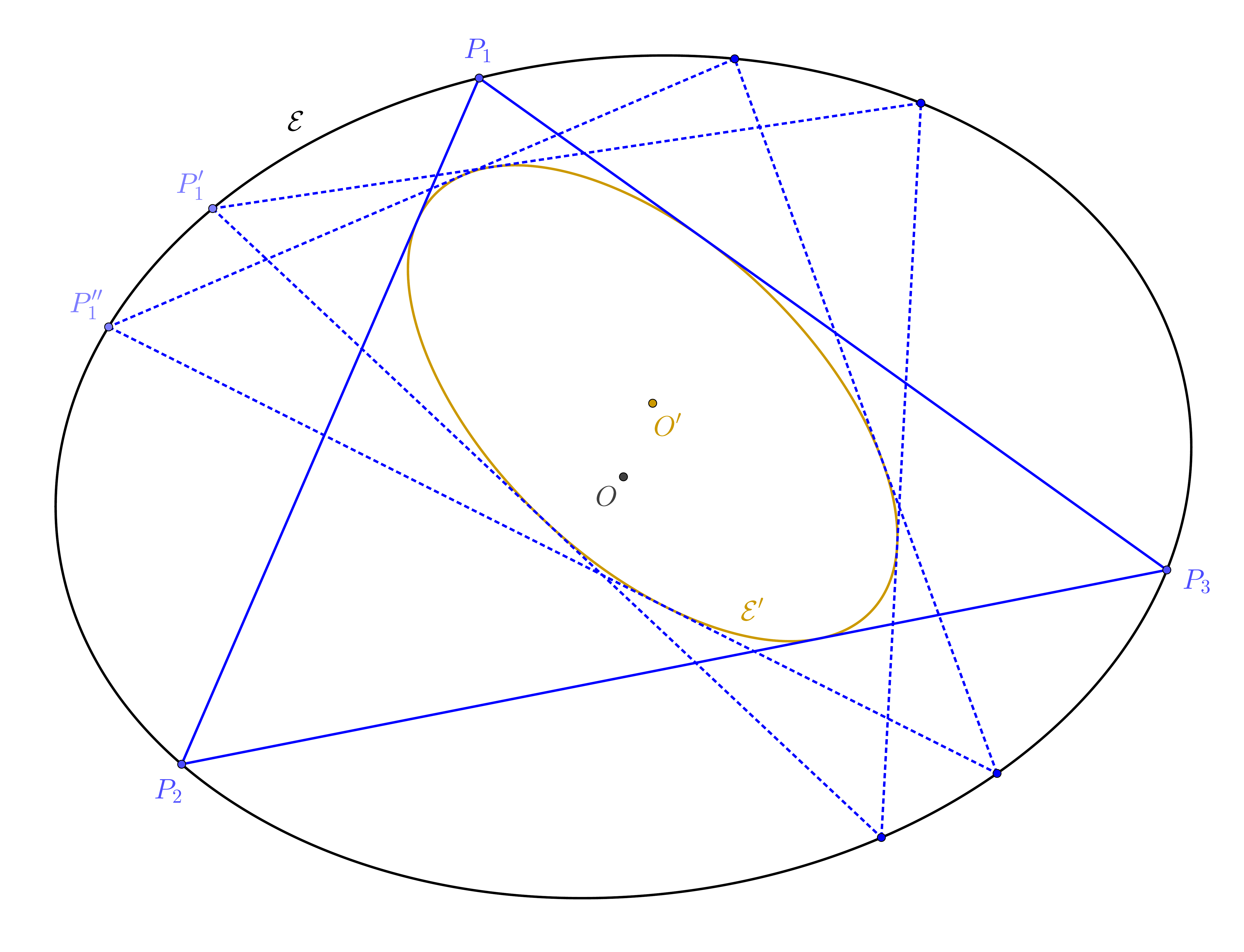}
     \end{subfigure}
     \rulesep
     \begin{subfigure}[c]{0.49\textwidth}
         \centering
         \includegraphics[width=\textwidth]{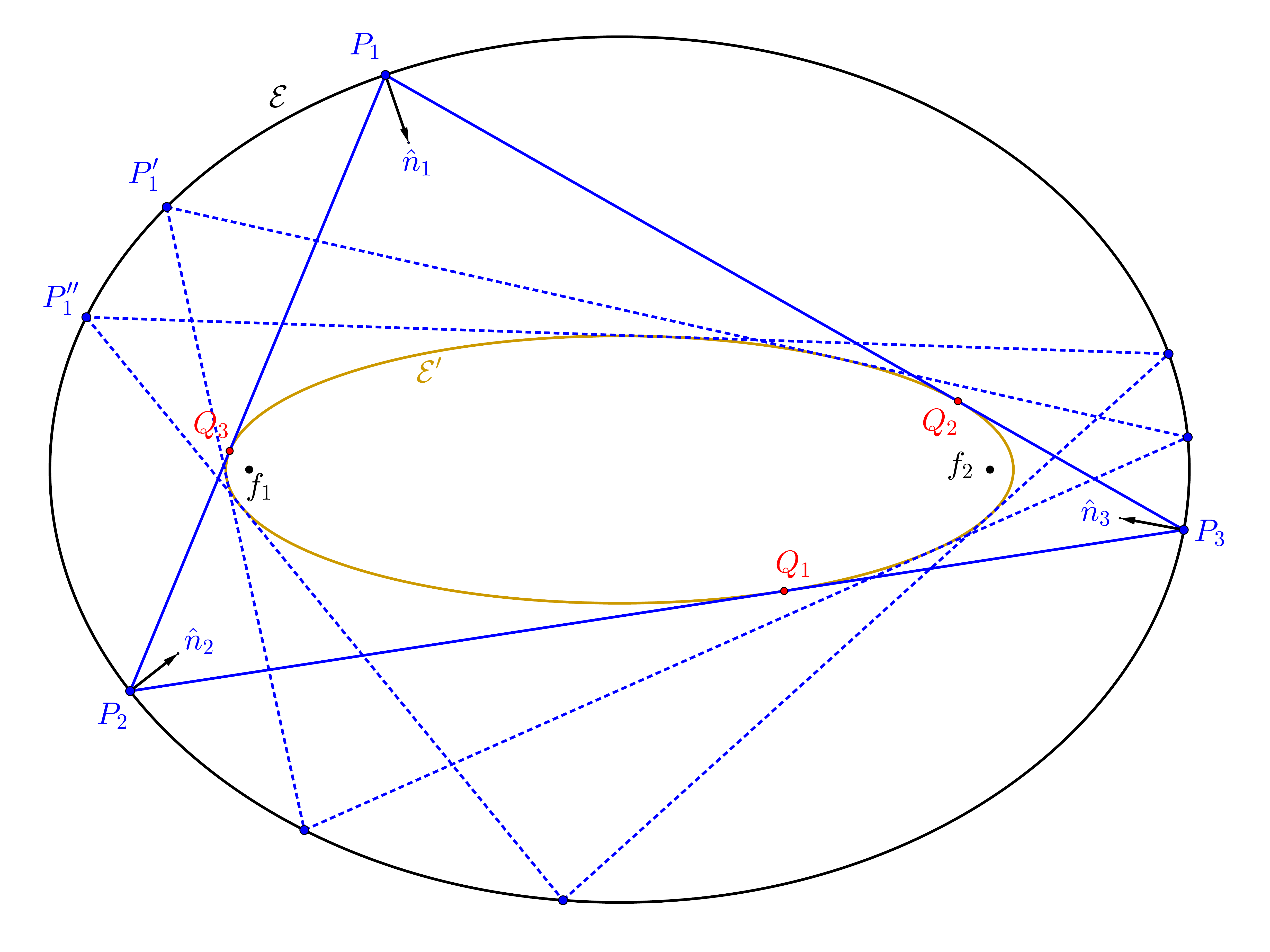}
    \end{subfigure}
    \caption{\textbf{Left:} Poncelet's closure theorem for triangles. \textbf{Right:} Ellipses $\E,\E'$ are confocal. Consecutive sides of the Poncelet family (blue) are bisected by the normals $\hat{n}_i$, and the perimeter is constant. Also shown are the three points of contact $Q_i$ with the inner ellipse, or caustic.}
    \label{fig:poncs}
\end{figure}

Depicted in \cref{fig:poncs}(left) is Poncelet's closure theorem in the special case of triangles. The theorem states that two conics\footnote{Recall these can be ellipses, hyperbolas, parabolas, and other degenerate specimens, see \cite[chapter 5]{gallier2011}.} $\E$ and $\E'$ are chosen so that a polygon can be drawn with all vertices on $\E$ and all sides tangent to $\E'$, then a {\em porism} of such polygons exists: any point $P$ on $\E$ can be used as an initial vertex for a polygon with identical incidence/tangency properties with respect to $\E,\E'$. For more details, see \cite{bos1987,centina15,dragovic11}.

Referring to \cref{fig:poncs}(right), a property-rich choice for $\E,\E'$ is when they are {\em confocal} ellipses, i.e., with shared foci. If such a pair admits a Poncelet porism\footnote{In general, finding such a pair requires that a certain ``Cayley'' determinant vanish, see \cite{dragovic11}.}, two immediate consequences ensue: (i) consecutive sides are bisected by the normal to $\E$ and Poncelet polygons can therefore be regarded as the periodic path of a particle bouncing elastically against $\E$ (this is known as the ``elliptic billiard'', see \cite{sergei91}), and (ii) all polygons in the porism have the same perimeter \cite{sergei91}. Dozens of other properties and invariants can be derived from these an interesting one being constant sum of the internal angle cosines, proved in \cite{akopyan2020-invariants,bialy2020-invariants}. For more properties of the confocal family, see \cite{garcia2020-new-properties,reznik2021-fifty-invariants}.

%%Poncelet f\cite{reznik2020-intelligencer}
%\cite{Hansmeyer2022}

\subsubsection*{Summary:} in \cref{sec:surfaces} we define a ruled surface based on Poncelet triangles, and discuss properties of its curvature. In \cref{sec:tangles} we study the link topology of space curves swept by  points of contact, and triangle centers. In \cref{sec:next-steps}, we list several unexplored experimental alternatives. To facilitate reproducibility, in \cref{app:jacobi,app:curvature} we include explicit expressions for both the Poncelet triangle parametrization and Gaussian and mean curvature. The pages listed in \cref{tab:resources} allow for live interaction with some objects mentioned herein.

{\small
\begin{table}
\centering
\begin{tabular}{|l|l|l|}
\hline
title & \texttt{http://observablehq.com/<.>} &  \\
\hline
PSTS Visualization (live) & \texttt{@esperanc/3-periodic-elliptical-billiards-3d-sweep} & \href{https://observablehq.com/@esperanc/3-periodic-elliptical-billiards-3d-sweep}{go} \\
PSTS Visualization (static) & \texttt{@dan-reznik/elliptic-billiard-triangle} & \href{https://observablehq.com/@dan-reznik/elliptic-billiard-triangle}{go} \\
Poncelet's Closure Theorem & \texttt{@dan-reznik/poncelet-iteration} & \href{https://observablehq.com/@dan-reznik/poncelet-iteration}{go}  \\
Jacobi's Elliptic Functions & \texttt{@dan-reznik/jacobi-elliptic-functions} & \href{https://observablehq.com/@dan-reznik/jacobi-elliptic-functions}{go} \\
\hline
\end{tabular}
\caption{Pages with interactive simulations to the various phenomena mentioned in the article.}
\label{tab:resources}
\end{table}
}

\section{A Poncelet Spatio-Temporal Surface (PSTS)}
\label{sec:surfaces}
\begin{figure}
     \centering
     \begin{subfigure}[c]{0.49\textwidth}
         \centering
         \includegraphics[width=\textwidth]{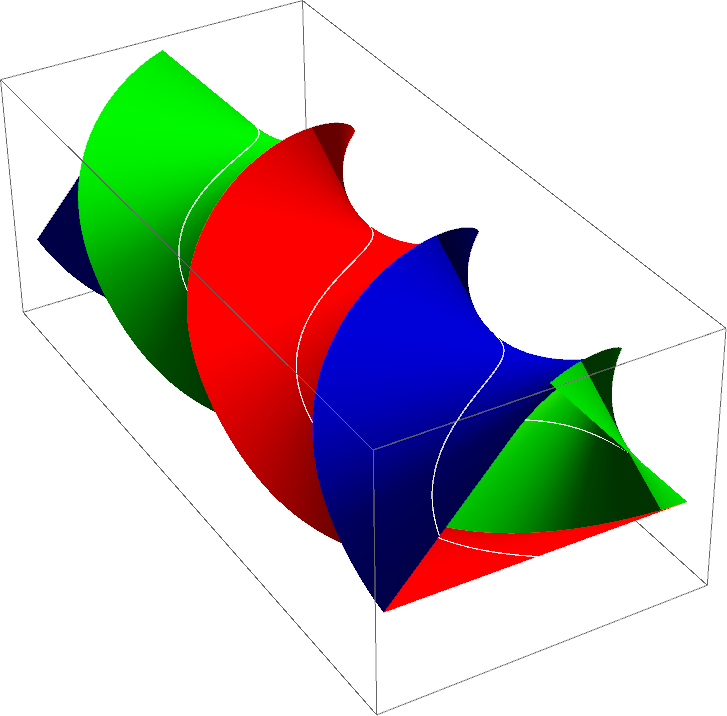}
     \end{subfigure}
     \rulesep
          \begin{subfigure}[c]{0.49\textwidth}
         \centering
         \includegraphics[width=\textwidth]{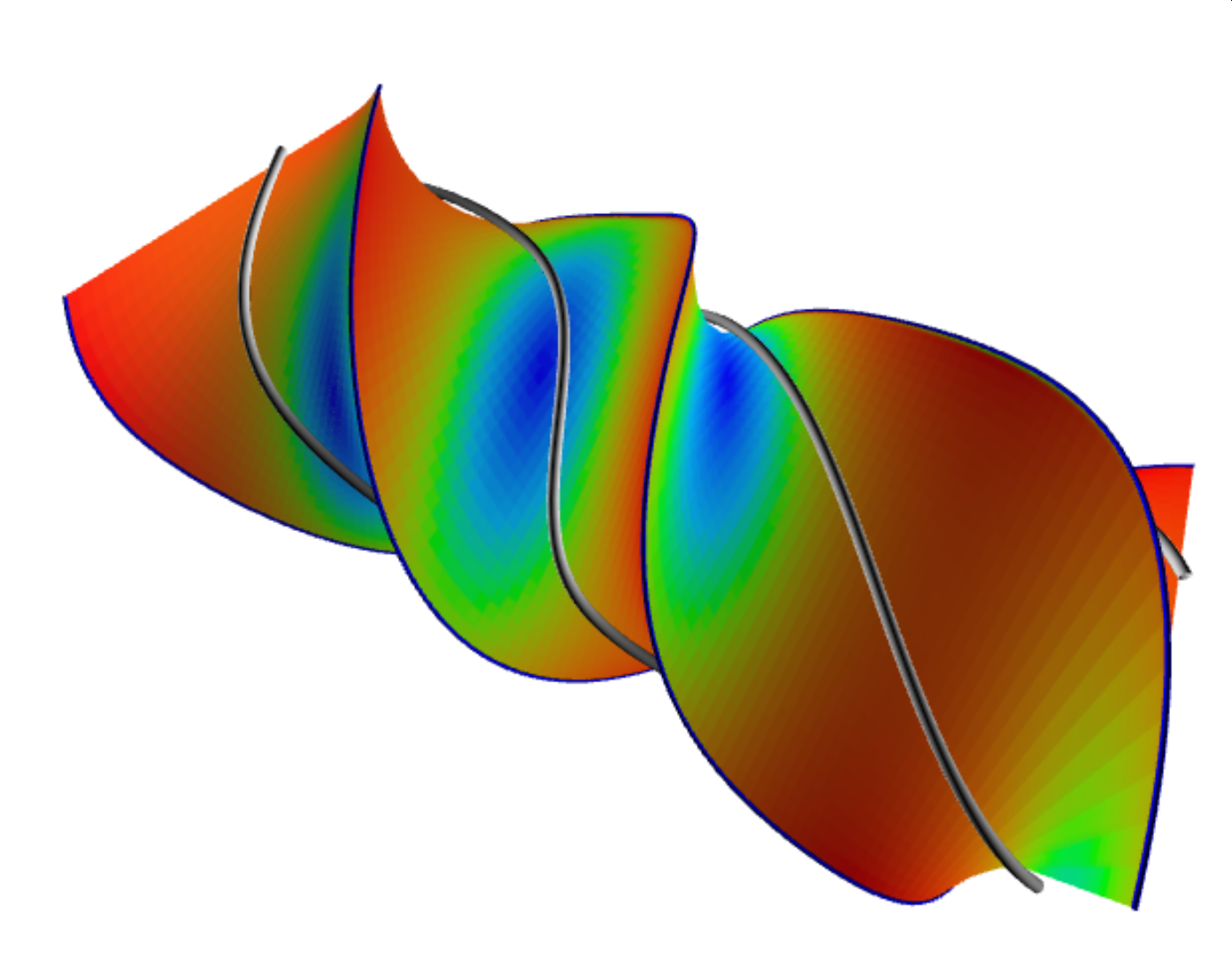}
     \end{subfigure}
\caption{\textbf{Left:} the PSTS swept by Poncelet triangles in the confocal family. It is a union of three ruled surfaces (red, green, and blue), each of which has negative curvature. Also shown is the 3d curve swept by the points of contact $Q_i$ (white). \textbf{Right:} the SPTS colored by Gaussian curvature. The center of the blue areas represent minima.}
\label{fig:corks}
\end{figure}

To achieve a homogeneous traversal of the Poncelet family, we parametrize it with Jacobi elliptic functions, as explained in \cref{app:jacobi}. Let $u$ be its parameter, $u\in[0,T]$, where $T$ is the period. Let $P_i(u)$ be a vertex of the family and $Q_i(u,v)$ be a point on edge $P_i(u) P_{i+1}(u)$, namely, $Q_i(u,v)=(1-v) P_i(u)+v P_{i+1}(u)$, $v\in [0,1]$. Referring to \cref{fig:corks}(left):

\begin{definition}
The Poncelet Spatio-Temporal Surface $\S$ (SPTS) is the union of the 3 parametric ruled surfaces $\S_i=[u,Q_i(u,v)]$, $i=1,2,3$. Note that the $u$ parameter is periodic.
\end{definition}

Recall that the Gaussian $\K$ (resp. mean $\H$) curvature) of a surface is the product (resp. average) of its principal curvatures, see \cite{spivak}. Referring to \cref{fig:corks}(right), and using the expressions in \cref{app:curvature}, we have derived rather long analytic expressions for both curvatures. Laborious analysis reveals that:

\begin{proposition} The three facets $\S_i$ of $\S$ are hyperbolic, i.e., each has negative Gaussian curvature everywhere.
\end{proposition}

\begin{figure}
\centering
     \begin{subfigure}[c]{0.49\textwidth}
         \centering
    \includegraphics[width=\textwidth]{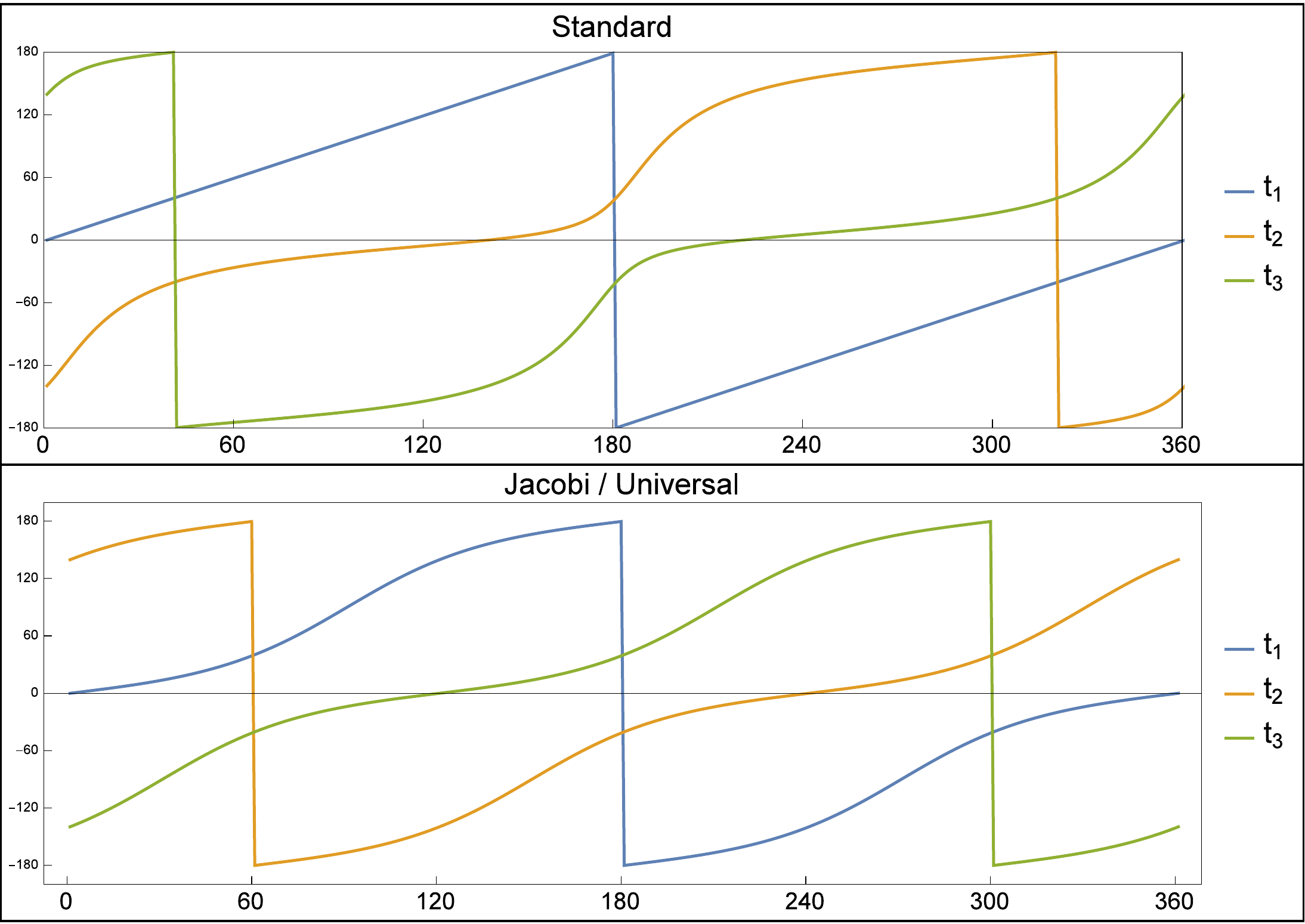}
\end{subfigure}  
   \hfill
   \begin{subfigure}[c]{0.49\textwidth}
    \centering
    \includegraphics[width=\textwidth]{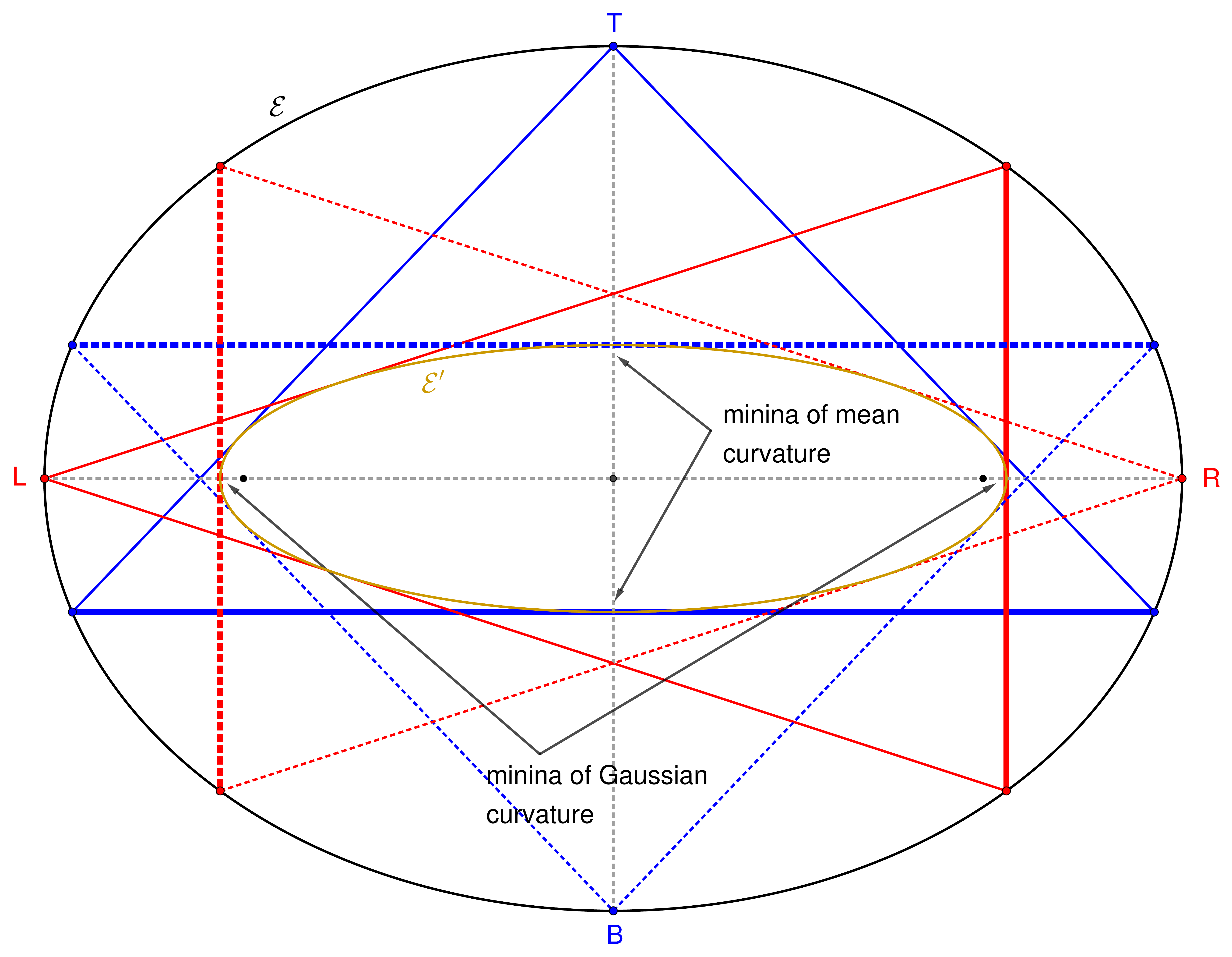}
   \end{subfigure}
\caption{\textbf{Left top:} under the ``standard parametrization'' $P_1(t)=[a\cos{t},b\sin{t}]$ the angular position of vertices of Poncelet triangles in the confocal pair are three different curves. \textbf{Left bottom:} under Jacobi's parametrization, the curves become 120-degree delayed copies of one another. \textbf{Right:} The confocal family has four isosceles triangles, with a vertex on either the top (T), bottom (B), left (L) or right (R) vertices of the outer ellipse $\E$. The Gaussian (resp. mean) curvature of the spatio-temporal surface have minima when its cross section is one of said isosceles triangles. The critical point occurs at the midpoint of the base (thick segment) when the apex is on L or R (resp. T or B).}
\label{fig:jacobi-isosceles}
\end{figure}

\begin{figure}
    \centering
    \includegraphics[width=\textwidth]{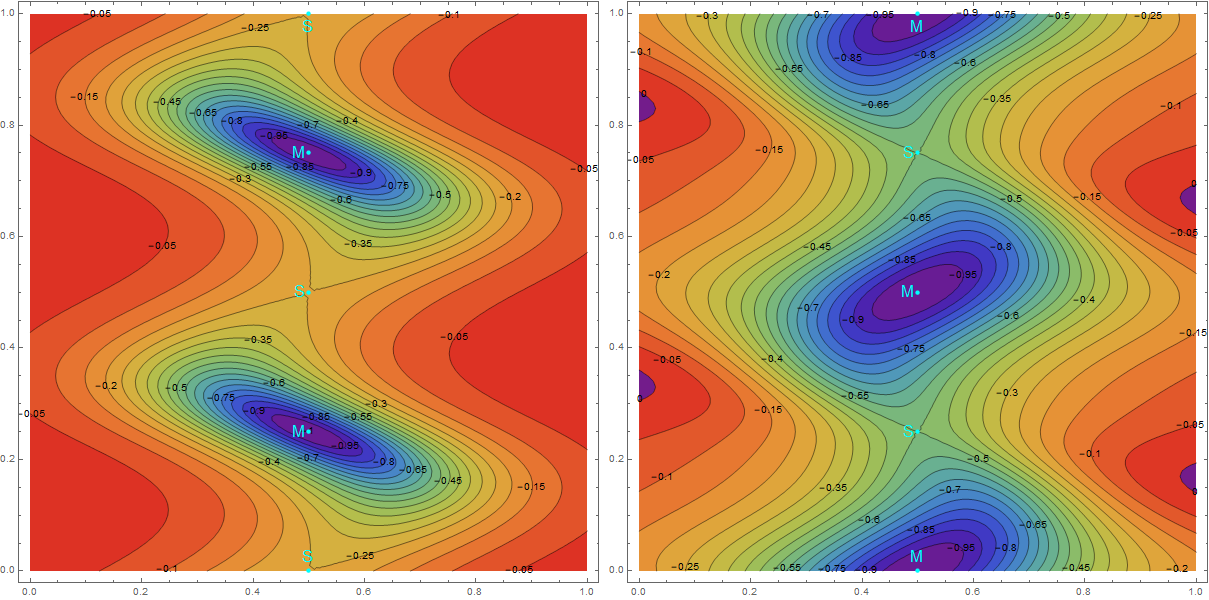}
    \caption{\textbf{Left:} Gaussian curvature of Jacobi-parametrized Poncelet triangles in the confocal family (horizontal is the position along a given side, and vertical is one revolution of the family. Points $M$ (resp. $S$) denote the curvature minima (saddle points). \textbf{Right:} The mean curvature, with $M,S$ as before.}
\label{fig:curvatures}
\end{figure}

Consider one facet $\S_1$ of $\S$. Let $Q_1=(0,1/2)$, $Q_2=(T/4,1/2)$, $Q_3=(T/2,1/2)$,  $Q_4=(3T/4,1/2)$. These four points correspond to the isosceles configurations shown in \cref{fig:jacobi-isosceles}(right). Referring to \cref{fig:curvatures}:

 \begin{proposition}
$Q_1$ and $Q_3$ (resp. $Q_2$ and $Q_4$) are non degenerate (Morse type) local minima (resp. saddlepoints) of $\mathcal{H}$. Conversely, $Q_2$ and $Q_4$ (resp. $Q_1$ and $Q_3$) are non degenerate (Morse type) local minima (resp. saddlepoints)  of $\mathcal{K}$.
\label{prop:HK}
\end{proposition}
 Analogous statements can be made for facets $S_2$ and $S_3$. It is worth noting that in general the critical points of Gaussian and Mean curvatures do not coincide. We currently think this is a feature of any Poncelet triangle family defined between a pair of concentric, axis-aligned ellipses.  

\begin{figure}
    \centering
    \includegraphics[width=\textwidth]{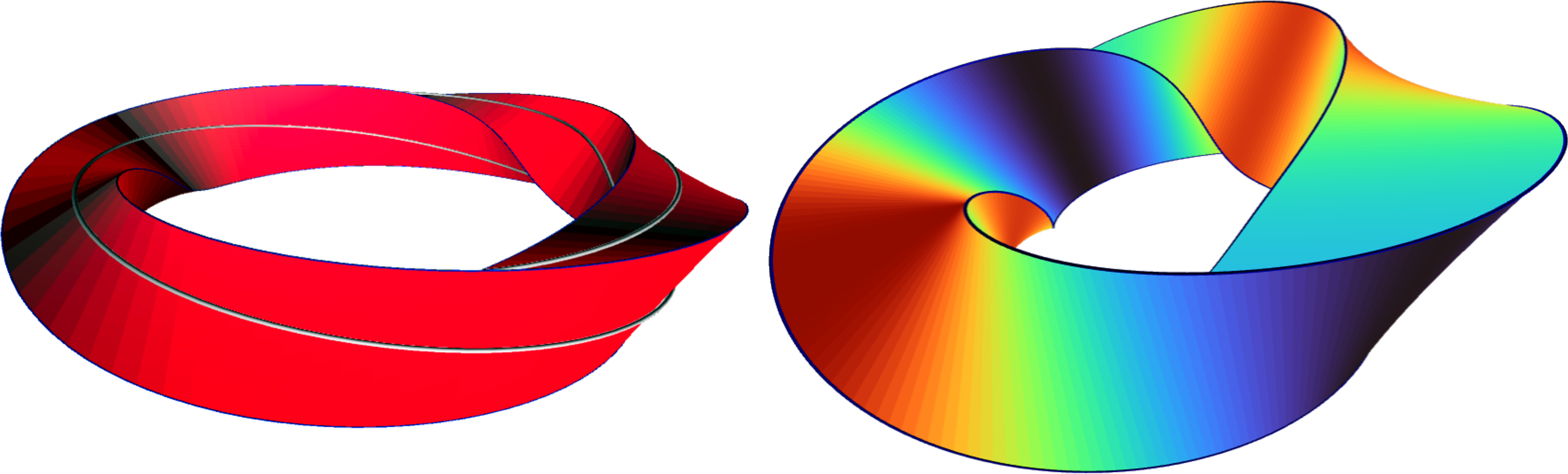}
    \caption{\textbf{Left:} identifying the $u=0$ with $u=T$ cross-sections of the PSTS, obtain an orientable Seifert surface \cite{wijk2006-seifert}. Also shown is the path of the contact points (white) with the caustic. \textbf{Right:} The same surface now colored by the torsion of straight lines elements sweeping the surface. Logo applications accepted.}
    \label{fig:seifert}
\end{figure}

\section{Space Curve Tangles}
\label{sec:tangles}
\begin{figure}
     \centering
     \begin{subfigure}[c]{0.49\textwidth}
     \centering
         \includegraphics[width=\textwidth]{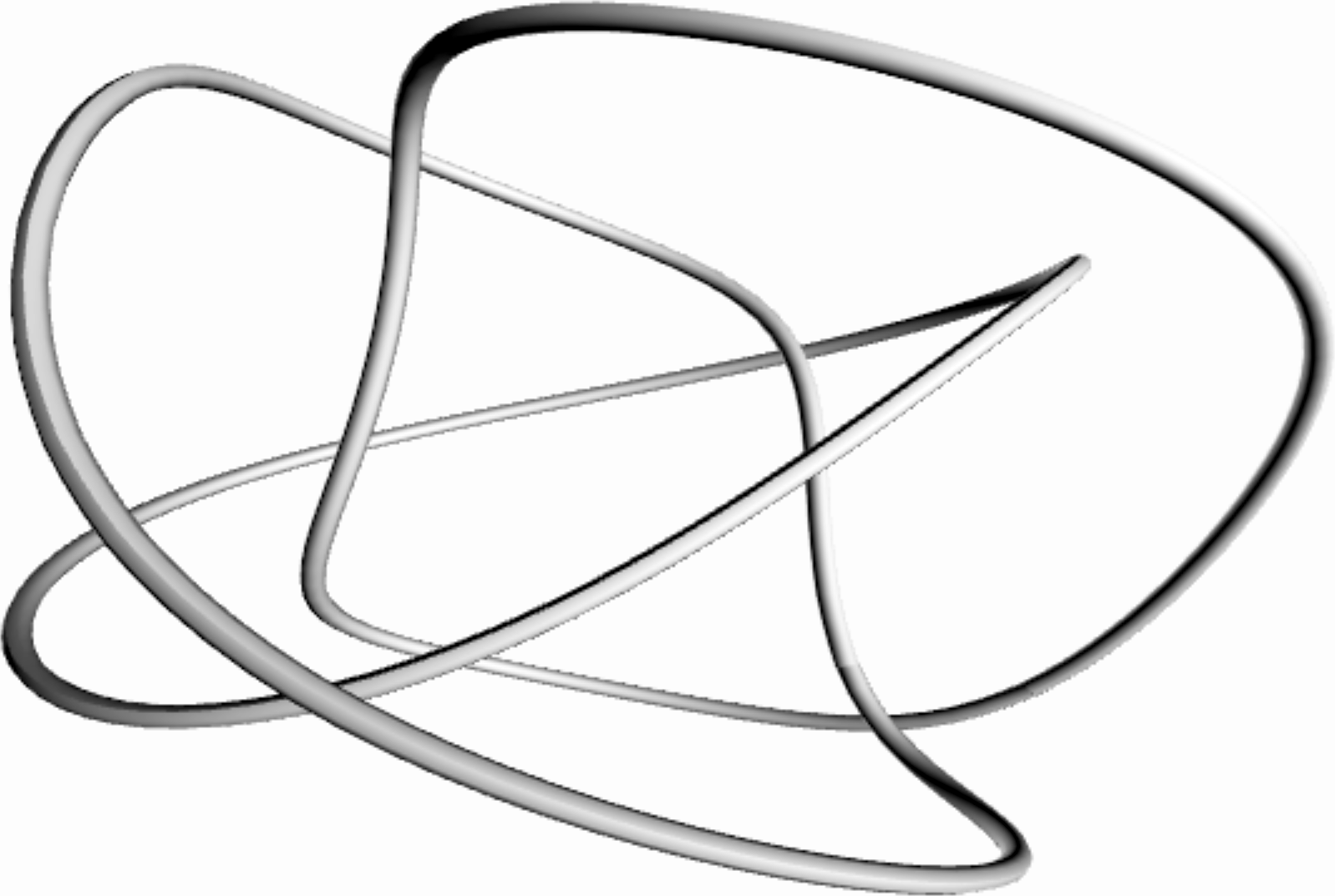}
                \end{subfigure}
                \rulesep
     \begin{subfigure}[c]{0.49\textwidth}
                \centering
         \includegraphics[width=\textwidth]{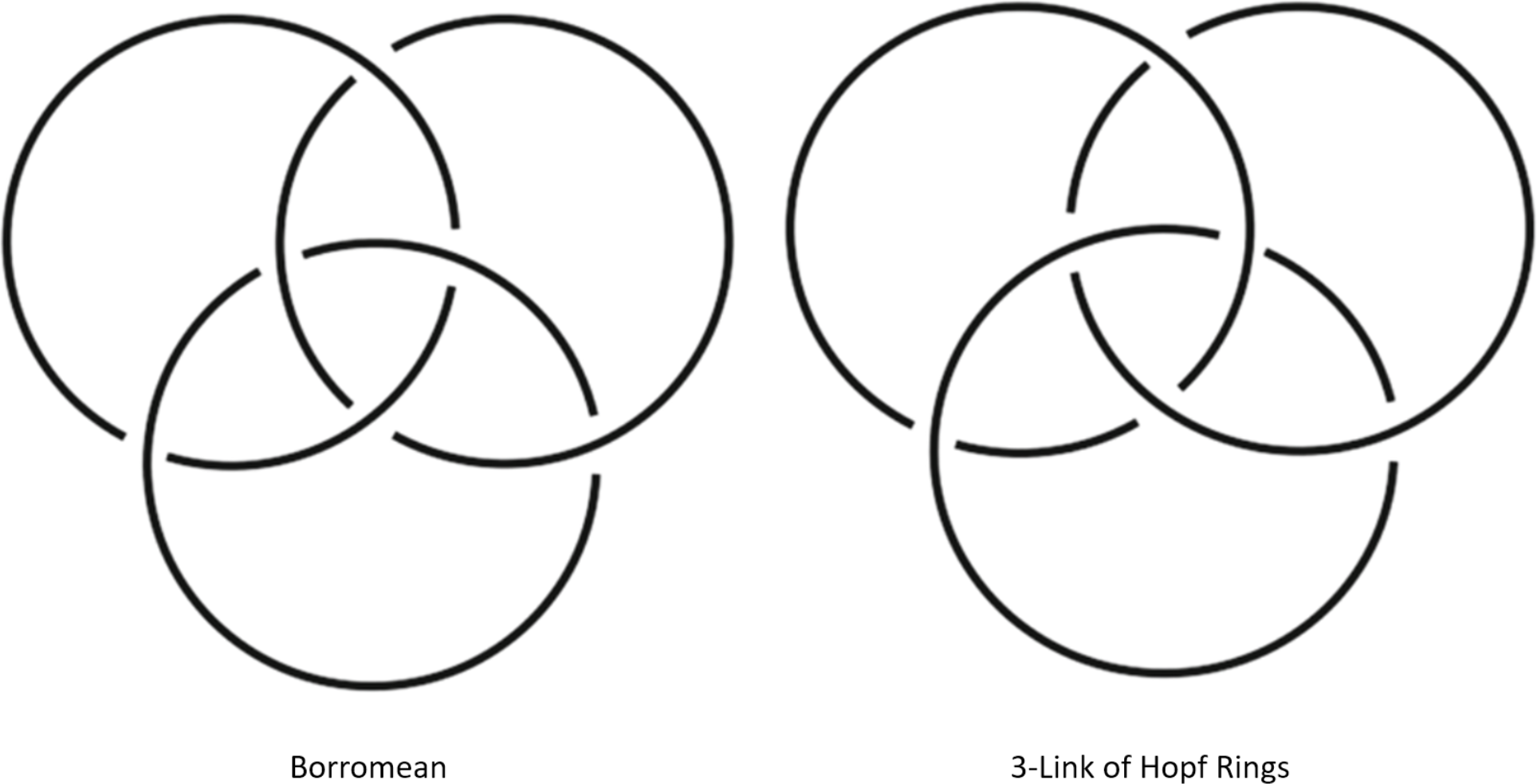}
         \end{subfigure}
        \caption{{Left:} The contact points of the identified PSTS sweep a triad of ``Hopf'' rings forming a 3-link tangle \cite{aravind1997-rings}. \textbf{Right:} Two type of 3-ring tangles: Borromean (left), and the Hopf 3-link (right) \cite{aravind1997-rings}, homeomorphic to the curves swept by the contact points. Note that by removing one of the rings in the former (resp. latter) case, the other two are free (resp. remain tangled). }
        \label{fig:tangles}
\end{figure}

Consider the surface obtained by identifying the $u=0$ and $u=T$ cross sections of $\S$, shown in \cref{fig:seifert}. Each contact point $Q_i$ of \cref{fig:poncs}(right) will sweep a wiggly ring; their union will form a {\em tangle} known as a 3-link of ``Hopf'' rings \cite{aravind1997-rings,rolfsen}, distinct from the Borromean tangle, see \cref{fig:tangles}. Indeed, the same tangle is swept by the 3 vertices of the family, and it is independent of the family being a Poncelet one. The surface whose boundary is a 3-link tangle is a type of {\em Seifert surface} \cite{wijk2006-seifert}.

%\begin{figure}
%    \centering
%    \includegraphics[width=.7\textwidth]{figures/fig_110_tangles_b.png}
%    \caption{Two types of 3-ring tangles: a Borromean tangle (left), and a 3-link (right) of ``Hopf'' rings \cite{aravind1997-rings}. Notice that in the former case, by removing one of the rings, the other two are disentangled, while this is not true in the latter case.}
%    \label{fig:tangles}
%\end{figure}

\begin{figure}
     \centering
       \includegraphics[width=\textwidth]{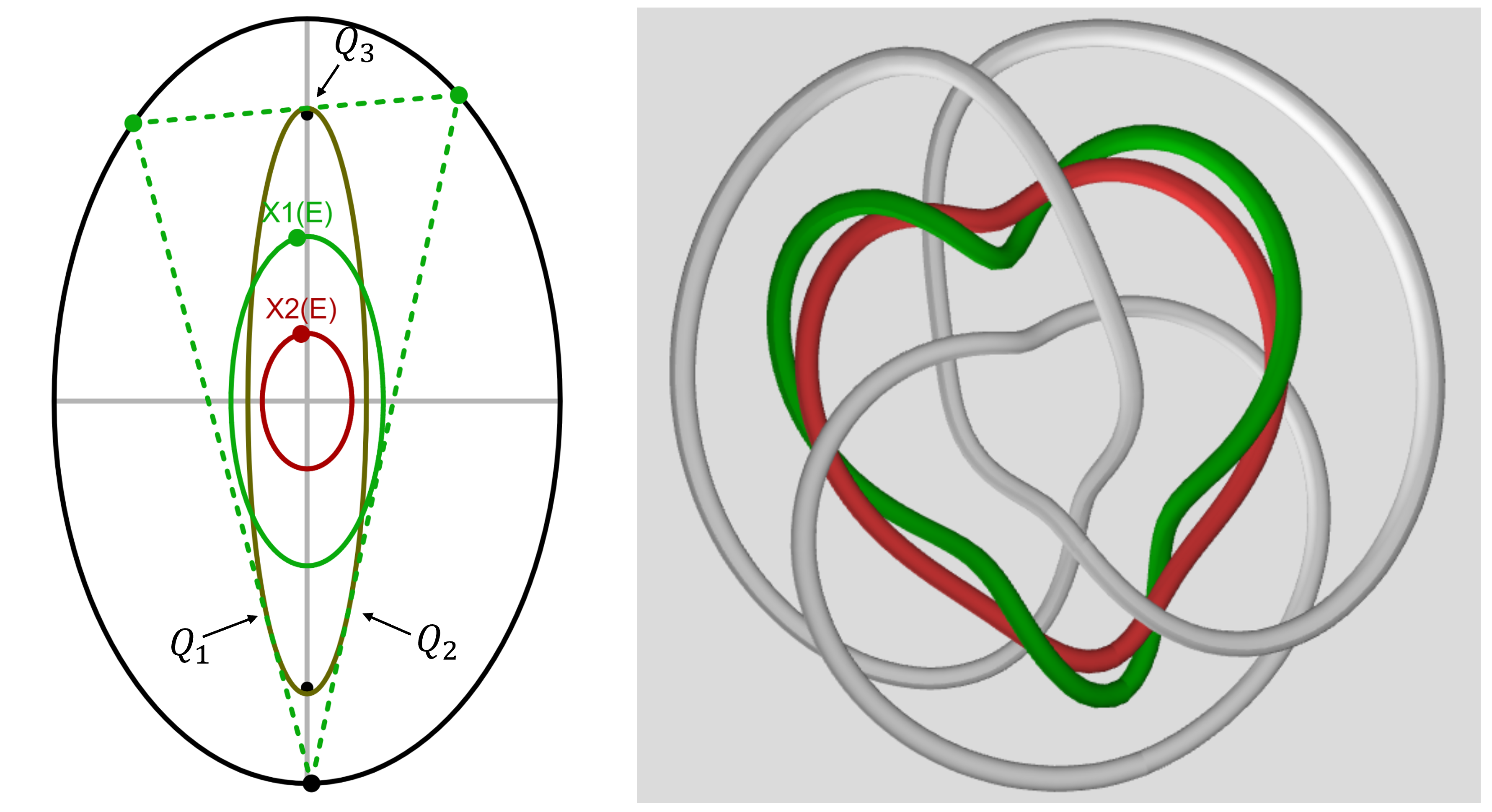}
       \caption{\textbf{Left:} the confocal family (rotated $90^o$), and the locus of the incenter $X_1$ (green) and barycenter $X_2$ (red). Also shown are the three contact points $Q_i$ with the caustic. \textbf{Right:} In the endpoint-identified PSTS, the space curves swept by $X_1$ (green) forms an individual a 2-link tangle with each individual contact point ring (gray). The same is true for the $X_2$ space curve (red). $X_1$ and $X_2$ form a link thrice twisted about each other.}
        \label{fig:tangles-xk}
\end{figure}

Referring to \cref{fig:tangles-xk}, more tangle topologies are obtained if one also considers the relative motion of notable points of the triangle (e.g., the incenter, the barycenter, etc.) over the Poncelet family. See \cite{reznik2020-intelligencer} for 2d analysis of such loci.

\section{Next Steps}
\label{sec:next-steps}
\begin{figure}
    \centering
    \includegraphics[width=\textwidth]{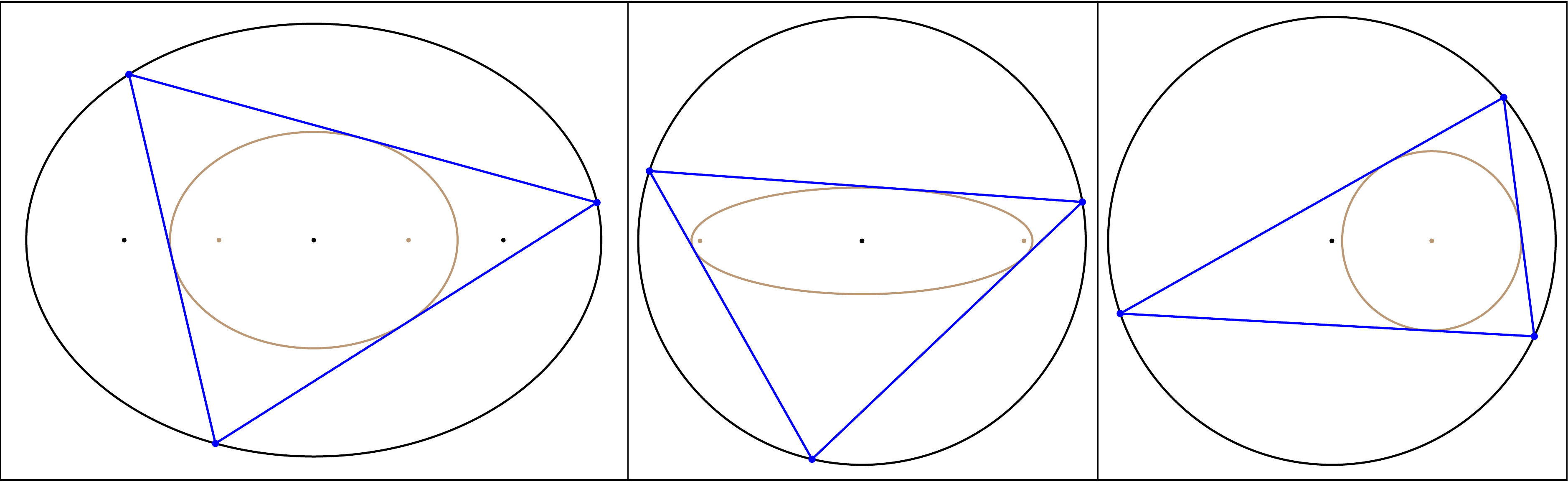}
    \caption{From left to right, three additional examples of Poncelet triangle families in (i) a homothetic pair of ellipses, (ii) inscribed in a circle and circumscribing a concentric ellipse, and (iii) interscribed between two non-concentric circles (aka., the ``bicentric'' pair). \href{https://youtu.be/14TQ5WlZxUw}{Video}}
    \label{fig:three-families}
\end{figure}

\begin{figure}
    \centering
    \includegraphics[ width=0.8\textwidth]{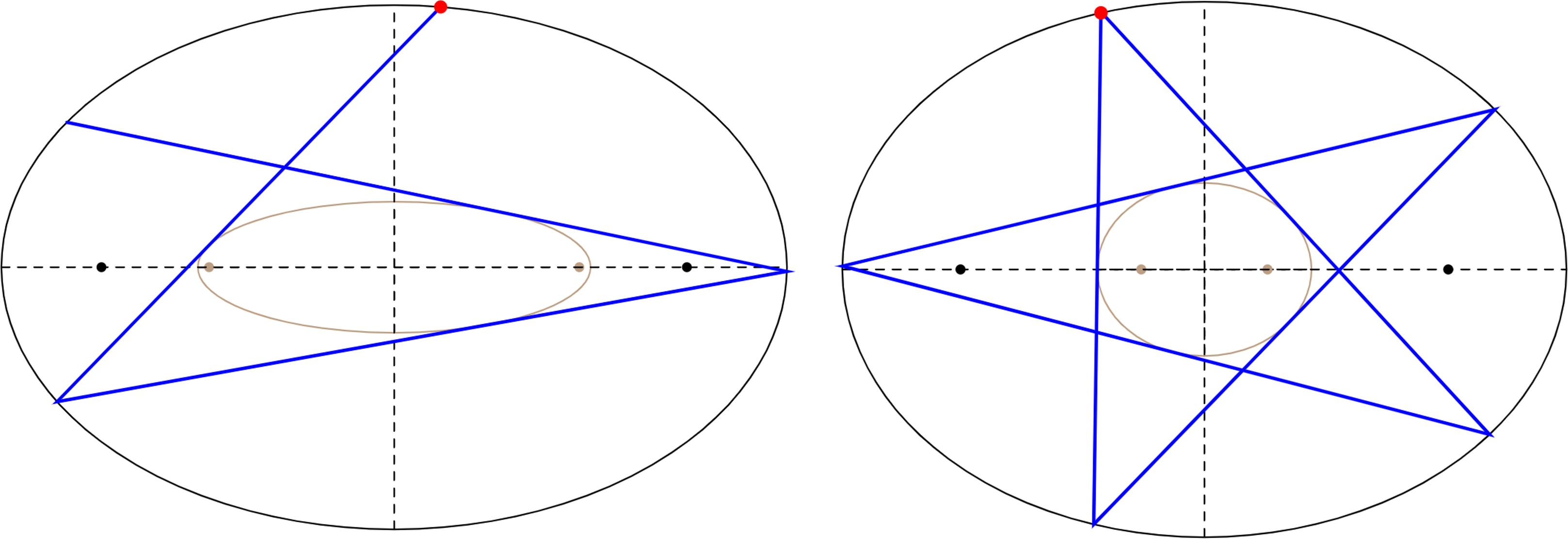}
    \caption{\textbf{Left:} a non-closing Poncelet 3-polyline. \textbf{Right:} a self-intersected $N=5$ Poncelet family (pentagrams).}
    \label{fig:open-pent}
\end{figure}

To continue this exploration one could consider:

\begin{compactenum}
    \item Different Poncelet families, see examples in \cref{fig:three-families};
    \item Picking a hyperbola or parabola for either $\E,\E'$, e.g., as in this \href{https://youtu.be/C8W2e6ftfOw}{video};
    \item Non-closing Poncelet polylines \cref{fig:open-pent}(left);
    \item Poncelet $N$-gons, $N>3$, including self-intersected ones as in \cref{fig:open-pent}(right). See \cite{garcia2020-self-intersected}.
    \item Families of derived triangles, e.g., the excentral, orthic, medial triangles \cite{mw}
\end{compactenum}

%\begin{figure}
%     \centering
%     \begin{subfigure}[c]{0.49\textwidth}
%         \centering
%         \includegraphics[trim=25 0 50 0,clip,width=\textwidth]{figures/fig_072_open_pentagram.png}
 %    \end{subfigure}
 %    \rulesep
 %    \begin{subfigure}[c]{0.49\textwidth}
 %        \centering
 %        \includegraphics[width=\textwidth]{figures/fig_075_excentral.png}
 %               \end{subfigure}
%\caption{\textbf{Left:} a non-closing Poncelet polyline (top) and a self-intersected $N=5$ (pentagram) family. \textbf{Right:} excentral triangles (green) to a reference Poncelet family (blue), in this case in a pair of confocal ellipses. In such a case, the excentral triangles will also be inscribed in an ellipse \cite{garcia2020-new-properties}.}
%\label{fig:extensions}
%\end{figure}

%\section*{Outline}
%\label{sec:outline}
%\input{900_outline}
 
%\section*{Acknowledgements}
%to do.

\appendix
\section{Jacobi Parametrization}
\label{app:jacobi}
We parametrize Poncelet triangles using Jacobi elliptic functions since an application of the Poncelet map corresponds to a unit translations in the argument of said functions. As seen in \cref{fig:jacobi-isosceles}(left), this entails that the angular position of vertices are identical, time-delayed functions.

Following the notation in \cite{armitage2006}, $k\in[0,1]$ denote the elliptic {\em modulus}\footnote{Mathematica (resp. Maple) expects $m=k^2$ (resp. $k$) for the second parameter to its elliptic functions.}:

\begin{definition}
The incomplete elliptic integral of the first kind $K(\varphi,k)$ is given by:
\begin{equation}
K(\varphi,k)=\int_0^{\varphi}\frac{d\theta}{\sqrt{1-k^2 \sin^2\theta}}
\label{eqn:02-ellipticK}
\end{equation}
%\end{definition}

%\begin{definition}
\noindent The {\em complete} elliptic integral of the first kind $K(k)$ is simply $K(\pi/2,k)$.
\end{definition}

\begin{definition}
The elliptic sine $\text{sn}$, cosine $\text{cn}$, and delta-amplitude $\text{dn}$ are given by:
\[
\text{sn}(u,k)=\sin\varphi,\;\;
\text{cn}(u,k)=\cos\varphi,\;\; 
\text{dn}(u,k)=\sqrt{1- k^2\sin^2\varphi}\]
Where $\varphi=\text{am}(u,k)$ as the {\em amplitude}, i.e., the upper-limit in the integral in \cref{eqn:02-ellipticK} such that $K(\varphi,k)=u$.
\end{definition} 

%A review of these functions appears in \cref{app:appD-jacobi-functions}.

%With this terminology the  Jacobi elliptic functions are defined    as follows:
%\begin{align*}
%K(\varphi,m)&=\int_0^{\varphi} \frac{dy}{\sqrt{1-m\sin^2y}}=u,\;\; \varphi=\text{am}(u,m)\\
%K(\varphi,k)=u\\ %
%\sn(u,m)&=\sin\varphi,\; \cn(u,m)=\cos\varphi,\; dn(u,m)=\sqrt{1-m\sin^2\varphi}%, \end{align*}

%As a note the reader, Maple accepts $k$ while Mathematica requires $m=k^2$ as the argument to elliptic functions.

%For example, in Maple we have: \[\sn(2,0.4)=\texttt{JacobiSN(2,.4)}=0.94569756\ldots\]
%and in Mathematica:

%\[ \sn(2,0.16)=\texttt{JacobiSN[2,.16]}=0.945698 \ldots \]
%and
%\[ \texttt{JacobiSN[2,.4]}=0.985090\ldots \]

\noindent As derived in  \cite{stachel2021-billiards}:
 
\begin{theorem}
A billiard N-periodic trajectory $P_i$ $(i=1,\ldots, N) $ of period $N$  with turning number $\tau$, where $\mathrm{gcd}(N,\tau) =1$ can be parametrized on $u$ with period $4K$ where:
\[ P_i=
%=\left[a\; \Jsn \left(u+\frac{4n\tau K}{N}, \frac{c}{a}\right), b\; \Jcn \left(u+\frac{4n\tau K}{N}, \frac{c}{a}\right)\right]\\
\left[-a\,\sn  \left(u+ i \Delta{u},  {m} \right) , b\,\cn  \left(u + i \Delta{u}, { m} \right)\right]
\]
\[ \text{with: } m=k^2=\frac{a_c^2-b_c^2}{a_c^2},\;\;\Delta{u}=\frac{4\tau K}{N},\;\;a= \sqrt{b^2+ a_c^2-b_c^2}, \;\; b=\frac{b_c}{\cn(\frac{\Delta{u}}{2}, m)}\]
\label{thm:stachel}
\end{theorem}

\section{Review: Gaussian and Mean Curvatures}
\label{app:curvature}
Let $\beta: M\to \mathbb R^3$ be a smooth immersion or embedding of a smooth oriented surface.
 The differential of $\beta$, $\beta_*$ is defined by
 $\beta_*(X)=d\beta(X)= X\beta$. The induced metric $g$, known as the {\em first fundamental form} is given by:
 $g(X,Y)=\langle \beta_*(X),\beta_*(Y)\rangle=\langle X\beta, Y\beta\rangle.$ Here $\langle,\rangle$ denotes the canonical inner product defining the Euclidean metric of $\mathbb R^3$.
 Consider an unit normal field $\mathcal N$ to the map $\beta$. The {\em second fundamental form} $\S:T_pM\to T_p M$ is defined by:
 \[X\mathcal{ N}=d{\mathcal N}(X) =-\beta_*(\S X)\]
 
% The structure equation is given by:
% $$
 %XY\beta=  \beta_*(\nabla_XY)+g(SX,Y){\mathcal N}.\;\;$$
 
 %The Codazzi  equation  is  given by:
 
 %\begin{equation}\aligned 
 %
 %(\nabla_X S)Y=& (\nabla_Y S)X,\;\;
 % \text{Codazzi equation}\\
 %K=& det(S)
 %
 %\endaligned
 %\end{equation}
 
 %Recall that
 %
 %$$D_X \beta_*(Y)=\beta_* (\nabla_XY)+\beta(X,Y)$$
 %where $\nabla_XY\in T_pM $ and %$\beta(X,Y)\in N_pM$.
 % D  is the usual derivative of $\mathbb R^3$,

 %Here $\nabla$ is the Levi-Civita connection of the metric $g$ and  $T_{\beta(p)}\mathbb R^3= \beta_*(T_pM)\oplus N_pM$.  
 
 The map $\S:T_pM\to T_pM$ is symmetric relative to the induced metric $g=\langle,\rangle_g$, i.e., 
 $\langle \S X,Y\rangle_g=\langle  X, \S Y\rangle_g$. The eigenvalues  $k_1\leq k_2$ of $\S$ are called the {\em principal curvatures} relative to $\mathcal{N}$ and the eigenspaces $e_i$ are called {\em principal directions}. The mean curvature $\mathcal{H}$ and Gaussian curvature $\mathcal{K}$ are given by \cite{spivak}:
 \[\mathcal{H}=(1/2)Tr(\S)=(k_1+k_2)/2,\;\;\; \mathcal{K}=det(\S)=k_1k_2\]

 %When $k_1<k_2$ the eigenspaces $e_1$ and $e_2$ are orthogonal relative to the metric $g$. The points where $k_1=k_2$ are called umbilic points and there all directions are principal directions.
 
% The integral curves of the principal directions are called lines of principal curvature or simply lines of curvature.
 
 In a local chart $(u,v)$ it follows that: 
 \[ \mathcal{H}=\frac{eG-2fF+Eg}{2(EG-F^2)}, \; \mathcal{K}=\frac{eg-f^2}{EG-F^2}\]
 where $I=Edu^2+2Fdudv+Gdv^2$ and $II=edu^2+2fdudv+gdv^2$
 are the first and second fundamental forms of the surface. Consider the Poncelet spatio-temporal surface $\S_1$. It follows that $\mathcal{H}=(H_n \Delta^{-\frac{3}{2}})/2$, and
 $\mathcal{K}=-(f/\Delta)^2$. Explicitly:
 {\small
 \begin{align*}
     H_n&=2   a   b   (d_4 + d_8)   (d_4   d_8 + s_4   s_8 - 1)   [(((a^2 - b^2)   s_4 - s_8   a^2)   d_4 + s_4   b^2   d_8)   (v - 1)   d_4 \\
     &- v   d_8   (-s_8   b^2   d_4 + (s_4   a^2 + (-a^2 + b^2)   s_8)   d_8)] + [    m^2 ab  (s_4 - s_8) \\
     &(2   d_4   s_4   s_8 + 2   d_4   s_8^2 - 2   d_8   s_4^2 - 2   d_8   s_4   s_8 - d_4 + d_8)   v \\
     &- a   b   (2  m^2  d_4     s_4^2   s_8 - 2   d_8   m^2   s_4^3 -  m^2  d_4    s_4 + m^2 d_8     s_4 - d_4   s_8 + d_8   s_4)]   [(a^2 - b^2)   s_4^2 - 2   s_4   s_8   a^2 \\
     &+ (a^2 - b^2)   s_8^2 - 2   b^2   (d_4   d_8 - 1)]\\
     \Delta &=[-2   [(m^2   s_4^2 + m^2   s_8^2 - 2   d_4   d_8 - 2)   v^2 + (-2   m^2   s_4^2 + 2   d_4   d_8 + 2)   v  
      + m^2   s_4^2 - 1]   [(s_8^2 - \frac{1}{2})   s_4^2 + s_8   (d_4   d_8 - 1)   s_4 -\frac{1}{2} s_8^2 \\
      &- d_8   d_4 + 1]
      b^2 
       - s_8^2 - 2   d_8   d_4 - s_4^2 + 2]   b^2  
      + b^2   (s_8 - s_4)^2 \\
     f&=-a   b   (d_4 + d_8)   (d_4   d_8 + s_4   s_8 - 1)
 \end{align*}
where $d_i=\dn(i K/3 + u, m)$, $s_i=\sn(i K/3 + u, m)$, $c_i=\cn(i K/3 + u, m)$, $i=4,8$, and $K$ is one quarter of the period in \cref{thm:stachel}, i.e., the complete elliptic integral of the first kind, see \cref{eqn:02-ellipticK}.
 }

%\section{Other Poncelet Triangle Families}
%\label{app:families}
%\input{220_cap_families}

\bibliographystyle{hacm}
\raggedright
\setlength{\baselineskip}{13pt}
\bibliography{999_refs,999_refs_rgk}
\end{document}